\documentclass[a4paper]{article}

\usepackage{hyperref}
\usepackage[vlined,boxruled,linesnumbered,algo2e,algosection]{algorithm2e}
\usepackage{lipsum}
\usepackage{amsmath,amsfonts}
\usepackage{cleveref}
\usepackage{graphicx}
\usepackage{epstopdf}
\ifpdf
  \DeclareGraphicsExtensions{.eps,.pdf,.png,.jpg}
\else
  \DeclareGraphicsExtensions{.eps}
\fi

\usepackage{pgfplotstable,filecontents}
\usepackage{pgfplots}
\usepackage{tikz}
\usetikzlibrary{fit,backgrounds,calc,pgfplots.statistics}

\usepackage{amsmath,amsthm}
\usepackage{amsopn}

\makeatletter
\newcommand*{\addFileDependency}[1]{
  \typeout{(#1)}
  \@addtofilelist{#1}
  \IfFileExists{#1}{}{\typeout{No file #1.}}
}
\makeatother

\makeatletter
\usetikzlibrary{decorations,backgrounds}
\def\pgfdecoratedcontourdistance{0pt}
\pgfset{
  decoration/contour distance/.code=%
    \pgfmathsetlengthmacro\pgfdecoratedcontourdistance{#1}}
\pgfdeclaredecoration{contour lineto closed}{start}{%
  \state{start}[
    next state=draw,
    width=0pt,
    persistent precomputation=\let\pgf@decorate@firstsegmentangle\pgfdecoratedangle]{%
    \pgfpathmoveto{\pgfpointlineattime{.5}
      {\pgfqpoint{0pt}{\pgfdecoratedcontourdistance}}
      {\pgfqpoint{\pgfdecoratedinputsegmentlength}{\pgfdecoratedcontourdistance}}}%
  }%
  \state{draw}[next state=draw, width=\pgfdecoratedinputsegmentlength]{%
    \ifpgf@decorate@is@closepath@%
      \pgfmathsetmacro\pgfdecoratedangletonextinputsegment{%
        -\pgfdecoratedangle+\pgf@decorate@firstsegmentangle}%
    \fi
    \pgfmathsetlengthmacro\pgf@decoration@contour@shorten{%
      -\pgfdecoratedcontourdistance*cot(-\pgfdecoratedangletonextinputsegment/2+90)}%
    \pgfpathlineto
      {\pgfpoint{\pgfdecoratedinputsegmentlength+\pgf@decoration@contour@shorten}
      {\pgfdecoratedcontourdistance}}%
    \ifpgf@decorate@is@closepath@%
      \pgfpathclose
    \fi
  }%
  \state{final}{}%
}
\makeatother
\tikzset{
  contour/.style={
    decoration={
      name=contour lineto closed,
      contour distance=#1
    },
    decorate}}

\newcommand{\SetAlgorithmStyle}{
  \setcounter{AlgoLine}{0}
  \SetKwData{Left}{left}\SetKwData{This}{this}\SetKwData{Up}{up}
  \SetKwInOut{Input}{Input}
  \SetKwInOut{Output}{Output}
  \ResetInOut{input}
  \SetKwComment{tcp}{//}{}
  \SetKwFor{For}{for}{}{end}
  \SetArgSty{}
  \DontPrintSemicolon
}
\ifpdf
\hypersetup{
  pdftitle={On the equivalence of the sparse Hermitian eigenvalue problem and hypergraph edge elimination},
  pdfauthor={K.~Kahl, B.~Lang}
}
\fi

\newtheorem{theorem}{Theorem}[section]
  \newtheorem{lemma}[theorem]{Lemma}

  \newtheorem{definition}[theorem]{Definition}
\newtheorem{remark}[theorem]{Remark}
\newtheorem{example}[theorem]{Example}

\newcommand{\OldStuff}[1]{{\color{red}#1}}
\renewcommand{\OldStuff}[1]{}

\newcommand{\DC}{divide-and-conquer}
\newcommand{\Hermit}{Hermit}
\newcommand{\RankOne}{rank-$1$}
\newcommand{\ROneMod}{\RankOne\ modification}
\newcommand{\ROneModified}{\RankOne-modified}

\newcommand{\AEE}{A_{E}}
\newcommand{\AVV}{A_{V}}
\newcommand{\Conj}[1]{\overline{#1}}
\newcommand{\Diag}{\operatorname{diag}}
\newcommand{\Dual}[1]{#1^{\ast}}
\newcommand{\ELower}{E^{\,\includegraphics[width=1.2ex]{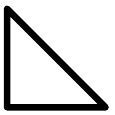}}}
\newcommand{\Herm}[1]{#1^{H}}
\newcommand{\IVE}{I_{V\!E}}
\newcommand{\IVEDual}{\Dual{\IVE}}
\newcommand{\MathC}{\mathbb{C}}
\newcommand{\MathR}{\mathbb{R}}
\newcommand{\Nnz}{\operatorname{nnz}}

\title{On the equivalence of the Hermitian eigenvalue problem and hypergraph edge elimination\thanks{
This work was partially funded by the DFG in the SFB TRR-55.}}

\newcommand{\OurAddr}{%
  University of Wuppertal,
  IMACM, Mathematics and Natural Sciences,
  D-42097 Wuppertal, Germany}
  
\author{%
  Karsten Kahl%
  \thanks{\OurAddr\ (\texttt{kkahl@math.uni-wuppertal.de}, \texttt{lang@math.uni-wuppertal.de})}
  \and
  Bruno Lang\footnotemark[2]%
  }
  
\begin{document}

\maketitle

\begin{abstract}
It is customary to identify sparse matrices with the corresponding adjacency or incidence graph. For the solution of linear systems of equations using Gaussian elimination, the representation by its adjacency graph allows a symbolic computation that can be used to predict memory footprints and enables the determination of near-optimal elimination orderings based on heuristics. The \Hermit{}ian eigenvalue problem on the other hand seems to evade such treatment at first glance due to its inherent iterative nature. In this paper we prove this assertion wrong by showing the equivalence of the \Hermit{}ian eigenvalue problem with a symbolic edge elimination procedure. A symbolic calculation based on the incidence graph of the matrix can be used in analogy to the symbolic phase of Gaussian elimination to develop heuristics which reduce memory footprint and computations. Yet, we also show that the question of an optimal elimination strategy remains NP-hard, in analogy to the linear systems case.
\end{abstract}




\section{Introduction}
\label{sec:intro}

The \DC\ algorithm is a well-known method for computing the eigensystem
(eigenvalues and, optionally, associated eigenvectors) of a \Hermit{}ian
tridiagonal matrix
  \cite{Cuppen1981,1997-Demmel,1987-DongarraSorensen%
  }.
It can be parallelized efficiently
  \cite{2011-AuckenthalerBlumEtAl,1987-DongarraSorensen%
        },
and even serially it is among the fastest algorithms available
  \cite{1999-AndersonBaiEtAl,1997-Demmel%
        }.

The method relies on the fact that if the eigensystem of a \Hermit{}ian
matrix $A_{0}$ is known, then the eigenvalues of a ``\ROneMod'' (or
``\RankOne\ perturbation'') of this matrix,
$A_{1} = A_{0} + \rho z \Herm{z}$, can be determined efficiently by
solving the so-called ``secular equation''
  \cite{BunchNielsenSorensen1978,1973-Golub%
        },
and $A_{1}$'s eigenvectors can also be obtained from those of $A_{0}$
  \cite{1994-GuEisenstat}.

In the tridiagonal case this can be used to zero out a pair of
off-diagonal entries $t_{k+1,k}$ and $t_{k,k+1} = \Conj{t_{k+1,k}}$ near
the middle of the tridiagonal matrix $T$ such that $T$ decomposes into
two half-size matrices and a \ROneMod,
\[
  T = \left[
        \begin{array}{@{}cc@{}}
          T_{1} &   0   \\
            0   & T_{2}
        \end{array}
      \right]
      \pm t_{k+1,k} z \Herm{z}
      ,
\]
with a vector $z$ containing nonzeros at positions $k$ and $k+1$ and
zeros elsewhere.
Having computed the eigensystems of $T_{1}$ and $T_{2}$ (by recursive
application of the same scheme), the eigensystem of $T$ is obtained from
these using the \RankOne\ machinery.

In this work we extend this method to a more general setting.
In
  \Cref{sec:see}
we show that the eigensystem of a \Hermit{}ian matrix can be computed
via a sequence of \ROneMod{}s, each of them removing entries
of the matrix until a diagonal matrix is reached.
  \Cref{sec:rankmod}
reviews some of the theory for \ROneMod{}s, as far as this is
essential for the following discussion.

While this approach in principle also works for full matrices, it
benefits heavily from sparsity.
In
  \Cref{sec:seenhg}
we show that the necessary work for a whole sequence of \ROneMod{}s
can be modelled in a graph setting, similarly to the fill-in arising in
direct solvers for \Hermit{}ian positive definite linear systems; cf.,
e.g.,
  \cite{2006-Davis,1981-GeorgeLiu%
        }.
However, the removal of nodes from the graph associated with the matrix
is not sufficient to fully describe the progress of the eigensolver;
here, the removal of \emph{edges} in hypergraphs
  \cite{Hypergraphs}
provides a natural description.

We present two ways to come back to node elimination.
In
  \Cref{sec:duality}
we consider the dual hypergraph, and in
  \Cref{sec:edgeedge}
we will see that the edge elimination is closely related to Gaussian
elimination for the so-called edge--edge adjacency matrix (and thus to
node elimination on the graph associated with that matrix).
In particular, an NP-completeness result will be derived from this
relation in
  \Cref{sec:NP}.
This result implies that, for a given sequence of \ROneMod{}s, it will
not be practical to determine an ordering of this sequence that is
optimal in a certain sense.

\OldStuff{%
  It is, however, not possible to determine the ``best'' sequence within
  a reasonable complexity; this will be discussed in \Cref{sec:NP}
  by showing that the search for optimal sequences is equivalent to
  the minimum-fill problem for suitable \Hermit{}ian linear systems,
  and the latter problem is known to be NP-complete
  \cite{1981-Yannakakis}.

  The paper is structured as follows. In~\Cref{sec:see} we first
  show that the \Hermit{}ian eigenvalue problem can be solved by
  successive edge elimination from the graph $G_{A}$ associated with
  $A$. In order to show that this approach is viable, we review results
  on \ROneModified\ eigenvalue problems in~\Cref{sec:rankmod}. In
  \Cref{sec:seenhg} we introduce the symbolic representation of the edge
  elimination procedure. To this end we introduce the concept of a
  hypergraph and develop the symbolic algorithm. In \Cref{sec:edgeedge}
  we discuss the duality of the developed process to the symbolic
  process of Gaussian elimination, which can be interpreted as
  successive node elimination in $G_{A}$.
}

Nevertheless, the hypergraph-based models allow to devise heuristics
for choosing among possible sequences of \ROneMod{}s such that the
overall consumption of resources is reduced.
In
  \Cref{sec:seecomp}
we discuss heuristics for the elimination orderings to reduce memory
footprint and computations.


Throughout the paper we assume that $A \in \MathC^{n \times n}$ is
\Hermit{}ian.
The presentation is aimed at sparse matrices, but ``sparsity'' is to be
understood in the widest sense, including full matrices.


\section{Successive edge elimination}%
\label{sec:see}

We first show that the \Hermit{}ian eigenvalue problem $AQ = Q\Lambda$
can be solved by a series of \ROneModified\ eigenvalue problems.
One way to do this is to have each \ROneMod\ remove one pair of nonzero
off-diagonal entries $a_{k,\ell}$ and $a_{\ell,k} = \Conj{a_{k,\ell}}$,
which in turn correspond to a pair of edges of the graph associated to
$A$.
Thus we first introduce the basic graph notation we require.

\begin{definition}
  The directed \emph{adjacency graph} $G_{A} = (V,E)$
  with vertex set $V$ and edge set $E$ that is associated with
  $A \in \MathC^{n \times n}$ is defined by
  \[
    V = \{1,\ldots,n\}
    \quad \text{and}\quad
    E = \{ (k,\ell) \in V^{2} \mid k \not= \ell, \ a_{k,\ell} \neq 0 \}
    .
  \]
\end{definition}

As our method treats matrix entries by conjugate pairs and maintains
\Hermit{}icity throughout, it is sufficient to consider only the lower
triangle of the matrix, corresponding to
$\ELower = \{ e = ( k, \ell ) \in E \mid k > \ell \}$.

\begin{definition}
  For each edge $(k,\ell) \in \ELower$ with
  $a_{k,\ell} = r_{k,\ell} \cdot e^{i\theta_{k,\ell}} \in \MathC$,
  where $r_{k,\ell} = | a_{k,\ell} |$,
  we define a \emph{vector representation $z_{(k,\ell)} \in \MathC^{n}$
  of the edge} by
  \[
    \left(
      z_{(k,\ell)}
    \right)_{j}
    =
    \begin{cases}
      1 & \text{if\ } j = \ell,\\
      e^{i\theta_{k,\ell}} & \text{if\ } j = k,\\
      0 & \text{else}.
    \end{cases}
  \]
\end{definition}

Using these vectors we can rewrite $A$ as a sum of \ROneMod{}s to a
diagonal matrix.

\begin{lemma}\label{lem:rankonerep}
  Let $A \in \MathC^{n\times n}$ be sparse and \Hermit{}ian with
  associated graph $G_{A}=(V,E)$.
  Then
  \begin{equation}\label{eq:rankonerep}
    A = D + \sum_{(k,\ell) \in \ELower}
              r_{k,\ell}\cdot z_{(k,\ell)} \Herm{z_{(k,\ell)}},
  \end{equation}
  where $D = \Diag(d_{1},\ldots,d_{n})$ with
  \begin{equation}\label{eq:di}
    d_{i}
    \ = \
    a_{i,i} - \sum_{( k, \ell ) \in \ELower,
                    k = i \; \mathrm{or} \; \ell = i}
                r_{k,\ell}
    \ = \
    a_{i,i} - \sum_{j=1,j\not=i}^{n} | a_{i,j} |
    .
  \end{equation}

  \begin{proof}
    For each edge $(k,\ell)$ with $k > \ell$, the \RankOne\ matrix
    $r_{k\ell} \cdot z_{(k,\ell)} \Herm{z_{(k,\ell)}}$ is nonzero only
    at the four positions $\{ \ell,k \} \times \{ \ell,k \}$,
    where we find
    \[
      \Big(
        r_{k,\ell} \cdot z_{(k,\ell)} \Herm{z_{(k,\ell)}}
      \Big)_{\{\ell,k\} \times \{\ell,k\}}
      =
      r_{k,\ell} \begin{bmatrix}
                   1 & e^{- i\theta_{k,\ell}}\\
                   e^{i\theta_{k,\ell}} & 1
                 \end{bmatrix}
      =
      \begin{bmatrix}
        r_{k,\ell} & \Conj{a_{k,\ell}}\\
        a_{k,\ell} & r_{k,\ell}
      \end{bmatrix} .
    \]
    Thus the $i$th diagonal entry is changed only by those edges
    starting or ending at node $i$, which gives the first equality
    in
      \cref{eq:di}.
    The second equality is a direct consequence of the definition
    of $\ELower$ and the \Hermit{}icity of $A$.
  \end{proof}
\end{lemma}

\begin{remark}
  The entries of $D$ in~\cref{eq:rankonerep} correspond to the lower
  bounds of the Gershgorin intervals.
  By defining $z_{(k,\ell)}$ differently one can also obtain a
  representation of $A$ similar to~\cref{eq:rankonerep} such that the
  entries of $D$ correspond to the upper bounds of the Gershgorin
  intervals.
\end{remark}

The solution of the \Hermit{}ian eigenvalue problem starting
from~\cref{eq:rankonerep} is now straight-forward.
Fixing an ordering of the edges, i.e., defining
$\ELower = \{ e_{1},\ldots,e_{|\ELower|}\}$ we have
\begin{equation}\label{eq:specialedgechosen}
  A = \left(
        D + r_{e_{1}} \cdot z_{e_{1}} \Herm{z_{e_{1}}}
      \right)
      + \sum_{j = 2}^{|\ELower|}
          r_{e_{j}} \cdot z_{e_{j}} \Herm{z_{e_{j}}}
  .
\end{equation}
Assuming that the eigendecomposition of the \Hermit{}ian matrix
$D + r_{e_{1}} \cdot z_{e_{1}} \Herm{z_{e_{1}}}$ has been computed,
\begin{equation*}
  D + r_{e_{1}} \cdot z_{e_{1}} \Herm{z_{e_{1}}} = Q_1 D_{1} \Herm{Q_1}
\end{equation*}
with $Q_{1}$ unitary, we can rewrite~\cref{eq:specialedgechosen} as
\begin{equation*}
  A = Q_1
      \Bigg(
        D_1 + \sum_{j = 2}^{|\ELower|}
                r_{e_{j}} \cdot
                  \left( \Herm{Q_{1}} z_{e_{j}} \right)
                  \Herm{\left( \Herm{Q_{1}} z_{e_{j}} \right)}
      \Bigg)
      \Herm{Q_{1}},
\end{equation*}
i.e., we eliminated edge $e_{1}$ from~\cref{eq:rankonerep}.
Successive elimination of the remaining $|\ELower|-1$ edges,
involving the vector
$\Herm{Q_{j-1}} \cdots \Herm{Q_{1}} \cdot z_{e_{j}}
 = \Herm{( \prod_{i=1}^{j-1} Q_{i} )} \cdot  z_{e_{j}}$ in step $j$,
finally yields the eigendecomposition of $A$,
\begin{equation*}
  A = \Bigg( \prod_{j=1}^{|\ELower|} Q_j \Bigg)
      \, D_{|\ELower|} \,
      \Herm{\Bigg( \prod_{j=1}^{|\ELower|} Q_j \Bigg)}.
\end{equation*}
This approach is summarized in Algorithm~\ref{alg:edgeElimination}.

\begin{algorithm2e}
  \SetAlgorithmStyle
  Write $A = D_{0} + \sum_{e\in \ELower}
                       r_{e} \cdot z_{e} \Herm{z_{e}}$\;
  Choose an ordering of the edges
    $e_{1},e_{2},\ldots,e_{|\ELower|}$\;
  Set $Q = I$\;
  \For{$j = 1$ \KwTo $|\ELower|$}{
    Calculate eigendecomposition of $D_{j-1} + r_{e_{j}} \cdot z_{e_{j}} \Herm{z_{e_{j}}}
               = Q_{j} D_{j} \Herm{Q_{j}}$\;
    \For{$i = j+1$ \KwTo $|\ELower|$}{
      $z_{e_{i}} = \Herm{Q_{j}} \cdot z_{e_{i}}$\;
    }
    $Q = Q Q_{j}$\;
  }
  \caption{Successive edge elimination.}
  \label{alg:edgeElimination}
\end{algorithm2e}

In order to be able to compute the eigendecomposition in this way we
need an efficient way to solve eigenproblems of the kind ``diagonal
plus \RankOne\ matrix.''
It is well known that these problems can be easily dealt with in terms of the
secular function, as we review in
  \Cref{sec:rankmod}.
In order to come up with a symbolic representation of the elimination
procedure we have to analyze the effect of the elimination of a
particular edge $e_{j}$ on the remaining edges.
This symbolic representation is developed in
  \Cref{sec:seenhg}.


\section{Computing eigenvalues of \RankOne-modified matrices}%
\label{sec:rankmod}

In order to clarify the main tool needed throughout the remainder of this
work we review some classical results about the eigenvalues of \RankOne\
perturbed matrices.
\OldStuff{%
  Let
  \begin{equation*}
    B = A + X,
  \end{equation*}
  where $A$ is \Hermit{}ian and $X$ is \Hermit{}ian positive
  semi-definite of rank $1$.}%
The results cited here date back to
  \cite{gantmaher1960oszillationsmatrizen}
and are also contained in
  \cite[pp.~94--98]{wilkinson-aep}.
They were later-on used in
  \cite{BunchNielsenSorensen1978,%
        Cuppen1981}
to formulate the \DC\ method for tridiagonal eigenproblems.

\begin{theorem}[{\cite[Theorem 1]{BunchNielsenSorensen1978}}]%
\label{theo:rankonepertubation}
	  Let $D + \rho z \Herm{z} = Q \Lambda \Herm{Q}$ be the
  eigendecomposition of the \RankOne-modified matrix, where
  $D = \Diag(d_{1},\ldots,d_{n}) \in \MathR^{n\times n}$
  with $d_{1} \leq d_{2} \leq \ldots \leq d_{n}$,
  $\|{z}\| = 1$, and $\rho > 0$.
  Then the diagonal entries of
  $\Lambda = \Diag(\lambda_{1},\ldots,\lambda_{n})$
  are the roots of the ``secular equation''
  \begin{equation}\label{eq:seculareq}
    f(\lambda) = 1 + \rho \sum_{j=1}^{n}
                            \frac{|z_{j}|^{2}}{d_{j} - \lambda}.
  \end{equation}
  More specifically, let the $\lambda_{j}$ be ordered,
  $\lambda_{1} \leq \lambda_{2} \leq \ldots \leq \lambda_{n}$.
  Then it holds
  \begin{equation}\label{eq:totalshift}
    \lambda_{j} = d_{j} + \rho \mu_{j}
      \quad \text{with} \quad
        0 \leq \mu_{j} \leq 1\text{\ for\ } j = 1,\ldots, n
        \quad\text{and}\quad
        \sum_{j = 1}^{n} \mu_{j} = 1.
  \end{equation}
\end{theorem}

There are two important consequences of Theorem~\ref{theo:rankonepertubation} found
in~\cite[pp.~94--98]{wilkinson-aep}.

\begin{lemma}\label{theo:DnCdeflation}
  Using the same notation as in Theorem~\ref{theo:rankonepertubation}
  we obtain the following.
  \begin{enumerate}
    \item In case the eigenvalues of $D$ are pairwise different we find
      that $\lambda_{j} = d_{j}$ if and only if $z_{j} = 0$. 
    \item In addition, if all $z_{j} \neq 0$, we find that
      $d_{j} < \lambda_{j} < d_{j+1}, j = 1,\ldots,n$
      ($d_{n+1} = \infty$).
    \item Assume there exists a multiple eigenvalue $d_{j}$ of $D$
      with multiplicity $k$; w.l.o.g.\ $d_{j-k+1}=\ldots=d_{j-1}=d_{j}$
      and $\|{z_{j-k+1,\ldots,j}}\| \neq 0$.
      Then we find
      \begin{equation*}
        \lambda_{i} = d_{i},\ i = j-k+1,\ldots,j-1,
        \quad \text{and} \quad
        d_{j} < \lambda_{j} < d_{j+1}
        \quad (d_{n+1} = \infty).
    \end{equation*}
  \end{enumerate}
\end{lemma}


  Lemma~\ref{theo:DnCdeflation} is one of the key algorithmic ingredients of the \DC\ algorithm for
tridiagonal eigenproblems and known in this context as ``deflation.''

As described in~\cite{1987-DongarraSorensen} %
and exploited in the implementation of the \DC\ method, the root-finding
problem of~\cref{eq:seculareq} is highly parallel and can be
efficiently solved by a modified Newton iteration using hyperbolae
instead of linear ansatz functions.

Recall that in our context the vector for the $j$th \ROneMod\
(elimination of $e_{j}$) is
$\Herm{( \prod_{i=1}^{j-1} Q_{i} )} \cdot z_{e_{j}}$.
Therefore,
  Lemma~\ref{theo:DnCdeflation}
implies that this elimination only requires the solution of the
secular equation in at most
\begin{equation}\label{eq:costEdgeElimination}
  N_{e_{j}} = \Nnz\Big( \Herm{(\prod_{i<j} Q_{i})} \cdot z_{e_{j}} \Big)
\end{equation}
intervals, where $\Nnz(v)$ is the number of nonzero entries of a
vector $v$.
That is, at most $N_{e_{j}}$ of the entries of $D_{j-1}$
(i.e., eigenvalue approximations) change from $D_{j-1}$ to $D_{j}$.
Further, by Theorem~\ref{theo:rankonepertubation}
we obtain that all eigenvalues move in the same direction and the
total displacement of these eigenvalues is given by
$r_{e_{j}} \cdot \|
                   \Herm{( \prod_{i < j} Q_{i} )} \cdot z_{e_{j}}
                 \|^{2}
 = 2 | a_{e_{j}} |$
because $r_{e_{j}} = | a_{e_{j}} |$ and the norm of the vector,
$\| z_{e_{j}} \| = \sqrt{2}$, does not change under the orthogonal
transformation $\Herm{(\prod_{i<j} Q_{i})}$.

Using the above reasoning, one would be able to estimate the cost
of the overall elimination process for a given ordering of the edges,
\emph{if} the number of nonzeros in the vectors
$\Herm{( \prod_{i<j} Q_{i} )} z_{e_{j}}$ could be \emph{predicted}.
In the following section we show how to do this.

Being able to analyze the influence of the ordering of the edges on
the complexity of the calculations (in terms of the number of roots of
the secular equations that need to be calculated) also allows us to
determine an ordering that leads to low overall cost.
This topic is discussed in
  \Cref{sec:NP,sec:seecomp}.


\section{Edge elimination, hypergraphs and
         edge elimination in hypergraphs}%
\label{sec:seenhg}

In
  \Cref{sec:see}
we have seen that the eigendecomposition of a \Hermit{}ian (sparse)
matrix $A$ can be obtained by successively eliminating the edges
$e_{1},e_{2},\ldots,e_{|\ELower|}$ of the graph $G_{A}$ associated
with the matrix $A$.

It is well known that in the context of Gaussian elimination for
\Hermit{}ian positive definite matrices, the effect of eliminating one
\emph{node} (corresponding to selecting a pivot row and doing the row
additions with this row) directly shows in the (undirected) graph
$G_{A}$: removing the node and connecting all its former neighbors
introduces exactly those edges that correspond to the new fill-in
produced by the row operations
  \cite{2006-Davis,1981-GeorgeLiu %
        }.
This allows to determine the nonzero patterns of the matrix during the
whole Gaussian elimination before doing any floating-point operation.

A similar thing can be done for the nonzero patterns of the vectors
$\Herm{( \prod_{i<j} Q_{i} )} z_{e_{j}}$ resulting from preceding
eliminations.
However, as we are eliminating \emph{edges}, the graph $G_{A}$ is
not adequate for this purpose.
We have to generalize the concept of a graph and use what is known in
the literature as a hypergraph~\cite{Hypergraphs}.

\begin{definition}
  An undirected \emph{hypergraph} $G = (V,E)$ is defined by
  a set of vertices $V = \{ v_{1},\ldots,v_{n}\}$ and
  a set of hyperedges $E = \{ e_{1}, \ldots, e_{m} \}$, where
  $\emptyset \not= e_{j} \subseteq V$.
\end{definition}

\begin{example}\label{ex:hypergraph} 
The hypergraph with vertex set $V = \{1,2,3,4,5\}$ and set of hyperedges 
$E = \{e_1,e_2,e_3,e_4\} = \left\{\{1,2,5\},\{2,3\},\{1,3,4,5\},\{3,4\}\right\}$ is depicted 
in Figure~\ref{fig:hypergraph_example}, where each edge in represented by 
a closed line that contains all its vertices.
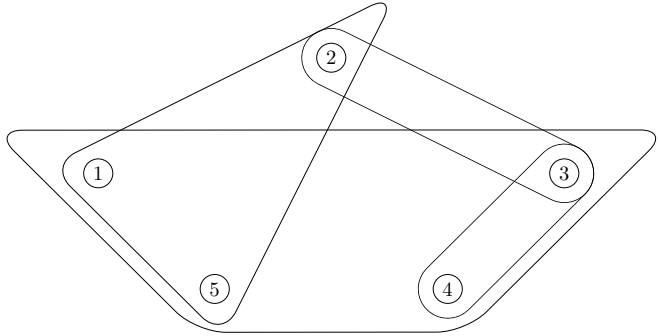
\begin{figure}
    \begin{center}
    \resizebox{.75\textwidth}{!}{
        \begin{tikzpicture}
            \path 
            (0,0) coordinate (V1) node[circle,inner sep=0pt,minimum width=.5cm,draw=black] {$1$}
            (4,2) coordinate (V2) node[circle,inner sep=0pt,minimum width=.5cm,draw=black] {$2$}
            (8,0) coordinate (V3) node[circle,inner sep=0pt,minimum width=.5cm,draw=black] {$3$}
            (6,-2) coordinate (V4) node[circle,inner sep=0pt,minimum width=.5cm,draw=black] {$4$}
            (2,-2) coordinate (V5) node[circle,inner sep=0pt,minimum width=.5cm,draw=black] {$5$};
            
            \path[preaction={contour=.5cm, rounded corners=.5cm,draw=black}] 
            (V1) -- (V2) -- (V5) -- cycle;
            \path[preaction={contour=.75cm, rounded corners=.5cm,draw=black}] 
            (V1) -- (V3) -- (V4) -- (V5) -- cycle;
            
            \pgfdeclarelayer{bg}    
            \pgfsetlayers{bg,main}

            \begin{pgfonlayer}{bg}
                \pgfsetblendmode{multiply}    

                \begin{scope}[transparency group]    
                \draw[double=none, double distance = 1cm, line join=round, line cap=round,fill=none] (V2) -- (V3) -- cycle;                
                \end{scope}

                \begin{scope}[transparency group]
                \draw[double=none, double distance = 1cm, line join=round, line cap=round,fill=none] (V3) -- (V4) -- cycle;
                \end{scope}
                
            \end{pgfonlayer}
        \end{tikzpicture}
    }
    \end{center}
    \caption{Drawing of the hypergraph defined in Example~\ref{ex:hypergraph} with $V = \{1,2,3,4,5\}$ and hyperedges $e_1 = \{1,2,5\}, e_2 = \{2,3\}, e_3 = \{1,3,4,5\}$ and $e_4 = \{3,4\}$.}\label{fig:hypergraph_example}
\end{figure}
\end{example}

\begin{remark}
  The possibility to have edges with more or less vertices than two is
  the only difference to the usual definition of an undirected graph.
  In particular, the graph $G_{A}$ can be considered as a hypergraph
  if we include each pair of edges $( k, \ell )$, $( \ell, k )$ only
  once, i.e., if we replace $E$ with $\ELower$.
\end{remark}

In order to analyze the nonzero pattern of the vector
$\Herm{( \prod_{i<j} Q_{i} )} z_{e_{j}}$ for the $j$th \ROneMod\ we
first note that this vector can be obtained in two ways:
``left-looking,'' when it is needed, by accumulating all previous
transformations $\Herm{Q_{i}}$ ($i < j$),
or ``right-looking,'' by applying each transformation $\Herm{Q_{i}}$,
once it has been computed, to all later $z_{j}$.
In the following discussion, as well as in Algorithm~\ref{alg:edgeElimination},
the right-looking approach is taken.

We now consider the effect of one such operation from the matrix/vector
point of view.
Let us assume that the edges are ordered and focus on the elimination of
the first edge, $e_{1}$.
Assume w.l.o.g.\ that $e_1 = \{1,2\}$.
By definition, $z_{e_{1}}$ has only two nonzero entries at the indices
$1$ and $2$, and thus due to Theorem~\ref{theo:rankonepertubation} and Lemma~\ref{theo:DnCdeflation}
we find
\begin{equation*}
  Q_{1} = \left[
            \begin{matrix}
              q_{11} & q_{12} & \\
              q_{21} & q_{22} & \\
                     &        & I_{(n-2) \times (n-2)}
            \end{matrix}
          \right].
\end{equation*}
Hence, for all edges $e_{j}$ with $e_{j} \cap e_{1} = \emptyset$ we have
$\Herm{Q_{1}} \cdot z_{e_{j}} = z_{e_{j}}$.
On the other hand, for all edges $e_{j}$ with
$e_{j} \cap e_{1} \neq \emptyset$ we find that
$\Herm{Q_{1}} \cdot z_{e_{j}}$ has
entries at the indices $e_{j} \cup e_{1}$.

The situation for the $i$th elimination step is similar.
Let the hyperedge $e_{j}$ denote the nonzero pattern, i.e., the set
of the positions of the nonzeros, of the current vector $z_{j}$ (after
the preceding transformations
$\Herm{Q_{i-1}} \cdots \Herm{Q_{1}} \cdot z_{j}$).
Then the transformed vector
$\Herm{Q_{i}} \cdot z_{j}$ will have nonzeros at the same positions
$e_{j}$ if $e_{j} \cap e_{i} = \emptyset$ and at positions
$e_{j} \cup e_{i}$ if the two hyperedges overlap.

\begin{remark}
  Strictly speaking this holds only if the transformation
  $\Herm{Q_{i}} \cdot z_{j}$ does not introduce new (``cancellation'')
  zeros in the vector.
  In the symbolic processing for sparse linear systems it is commonly
  assumed that this does not happen; we will do so as well.
\end{remark}

We summarize the above observation in the following theorem.

\begin{theorem}\label{def:HypergraphElimination}
  Let $G = (V,E)$ be an undirected hypergraph with $E \neq \emptyset$.
  Let $x \in E$ be the edge to be eliminated, and let
  \[
    E = E_{x} \cup E_{\not x},
    \quad \text{where} \quad
    \begin{cases}
      E_{x} = \{e\in E \mid e \cap x \neq \emptyset\} \text{\ and}\\
      E_{\not x} = \{e\in E \mid e\cap x = \emptyset\}.
    \end{cases}
  \]
  Then the hypergraph after elimination of $x$ is given by
  $\widetilde{G} = (\widetilde{V},\widetilde{E})$ with
  $\widetilde{V} = V$ and
  $\widetilde{E} = \{e\cup x, e\in E_{x}\setminus \{x\}\}
                   \cup E_{\not x}$.
\end{theorem}

Now it is easy to show that the subsequent elimination of all edges to
compute the eigendecomposition as described in
  \Cref{sec:see}
is equivalent to the elimination of all edges in the same ordering from
the (hyper)graph $G_A$ as defined here.
Thus it is natural to discuss questions such as complexity and optimal
edge orderings in the ``geometrical'' context of these graphs as it has
been successfully done for the solution of linear systems (e.g.,
optimal node orderings to reduce fill-in).

\begin{remark}
\label{rem:MoreGeneralSettings}
  In the above discussion we have assumed that each step of the
  algorithm eliminates a ``true edge'' $e = ( k, \ell )$, zeroing a pair of
  matrix entries $a_{k,\ell}$ and $a_{\ell,k}$.
  However this is not mandatory.
  Note that Theorem~\ref{def:HypergraphElimination}
  describes the evolution of the nonzero patterns also if the eliminated
  edge $x$ is a hyperedge as well, with more than just $2 \times 2$
  matrix entries being touched by the corresponding \ROneMod.
  In addition, the (off-diagonal) matrix entries at the positions
  $x \times x$ need not be zeroed out completely with the elimination.
  This allows for more general elimination strategies, including the
  extremes
  \begin{itemize}
    \item each \ROneMod\ zeroes one off-diagonal pair of matrix
      entries (cf.\ \Cref{sec:see}), and
    \item the $i$th \ROneMod\ zeroes the whole $i$th column and row
	of the matrix; this typically leads to the minimum
	\emph{number} of \ROneMod{}s, but according to the above the
      operations $\Herm{Q_{i}} \cdot e_{j}$ will make the vectors
	dense very quickly,
  \end{itemize}
  as well as many intermediate variants.
  For example, if the underlying model leads to low-rank off-diagonal
  blocks in the matrix then these can be removed with a reduced number
  of steps: for a size-($r \times s$) block of rank $\rho$, $\rho$
  \ROneMod{}s (with identical hyperedges) are sufficient instead of
  $r \cdot s$.
  We will come back to this generalization in
    \Cref{sec:NP}.
\end{remark}


\section{Duality between edge elimination and node elimination}%
\label{sec:duality}

%

In this section we will show that edge elimination can also be
expressed as \emph{node} elimination in a suitable graph.
This requires a few preparations.

\begin{definition}\label{def:incidence-matrix}
  Let $G = (V,E)$ be a hypergraph with
  nodes $V = \{v_{1},\ldots,v_{n}\}$ and
  hyperedges $E = \{e_{1},\ldots,e_m\}$.
  The \emph{(node--edge) incidence matrix}
  $\IVE \in \MathR^{|V| \times |E|}$
  of $G$ is then defined by 
  \[
    \left( \IVE \right)_{ij}
    =
    \begin{cases}
		1, & \text{if\ } v_{i} \in e_{j}, \\
      0, & \text{else},
    \end{cases}
  \]
  and the \emph{adjacency matrices of the hypergraph} are given by
  \[
    \begin{array}{r@{\;}c@{\;}l@{\quad}l}
      \AVV & = & \IVE \cdot \IVE^{T} \in \MathR^{|V| \times |V|}
           & \text{(\emph{vertex--vertex adjacency matrix})}, \\[1ex]
      \AEE & = & \IVE^{T} \cdot \IVE \in \MathR^{|E| \times |E|}
     & \text{(\emph{edge--edge adjacency matrix})}.
    \end{array}
 \]
\end{definition}

The latter two names are explained by the following lemma.

\begin{lemma}
  \label{lem:incidence-adjacency}
  Given a hypergraph $G = (V,E)$, its adjacency matrices have the
  properties
  \[
    \left( \AVV \right)_{ij}
    \neq 0 \quad \text{iff} \quad
    \text{there exists $e \in E$ such that $v_{i},v_{j} \in e$,}
  \]
  i.e., nodes $v_i$ and $v_{j}$ are connected by at least one hyperedge,
  and
  \[
    \left( \AEE \right)_{ij}
    \neq 0 \quad \text{iff} \quad
    \text{there exists $v \in V$ such that $v \in e_{i} \cap e_{j}$,}
  \]
  i.e., the hyperedges $e_{i}$ and $e_{j}$ share at least one node $v$.

  \begin{proof}
    Follows immediately from Definition~\ref{def:incidence-matrix}
    and the calculation of matrix--matrix products due to
    \[
      \left( \IVE \cdot \IVE^{T} \right)_{ij}
      =
      \sum_{k = 1}^{|E|}
        \left( \IVE \right)_{i,k} \cdot \left( \IVE^{T} \right)_{k,j}
      =
      \sum_{k=1}^{|E|}
        \left( \IVE \right)_{i,k} \cdot \left( \IVE \right)_{j,k}
      \; ;
    \]
    similarly for $\AEE$.
  \end{proof}
\end{lemma}

We also note that the transpose of the incidence matrix of $G$,
$\IVE^{T}$, is also the incidence matrix $\IVEDual$ of the dual of the
hypergraph, which is defined as follows.

\begin{definition}
  Let $G=(V,E)$ be a hypergraph with nodes $V=\{v_1,\dots,v_n\}$ and
  hyperedges $E=\{e_1,\dots,e_m\}$.
  Then the \emph{dual of $G$} is a hypergraph
  $\Dual{G} = ( \Dual{V}, \Dual{E} )$ with
  nodes $\Dual{V} = \{ \Dual{v_{1}}, \dots, \Dual{v_m} \}$ and
  hyperedges $\Dual{E} = \{ \Dual{e_{1}}, \dots, \Dual{e_n} \}$
  such that
  \[
    \Dual{e_i} = \{ \Dual{v_j} \in \Dual{V} \,:\, v_i \in e_j \}.
  \]
\end{definition}

By construction, edge elimination in a hypergraph is equivalent to
node elimination in its dual, as can be seen in the following small
example as well.

\begin{example}
    Consider the hypergraph $G=(V,E)$ of Example~\ref{ex:hypergraph_example} and its dual
  $\Dual{G} = (\Dual{V},\Dual{E})$ given by their incidence matrices
  $\IVE$ and $\IVEDual$, respectively:
  \[
    \IVE = \left[
             \begin{array}{cccc}
               1 & 0 & 1 & 0 \\
               1 & 1 & 0 & 0 \\
               0 & 1 & 1 & 1 \\
               0 & 0 & 1 & 1 \\
               1 & 0 & 1 & 0
             \end{array}
           \right]
    , \qquad 
     \IVEDual =
     \IVE^{T} = \left[
                  \begin{array}{ccccc}
                    1 & 1 & 0 & 0 & 1 \\
                    0 & 1 & 1 & 0 & 0 \\
                    1 & 0 & 1 & 1 & 1 \\
                    0 & 0 & 1 & 1 & 0
                  \end{array}
                \right]
    .
  \]
  If we eliminate edge $e_1$ in $G$ or, equivalently, node $\Dual{v_1}$
  in $\Dual{G}$, then the resulting hypergraphs $\widetilde{G}$ and
  $\widetilde{\Dual{G}}$ are given by
  \[
    \widetilde{\IVE}
    = \left[
        \begin{array}{cccc}
          & \mathbf{1} & 1 & 0 \\
          & 1 & \mathbf{1} & 0 \\
         \phantom{0} & 1 & 1 & 1 \\
          & 0 & 1 & 1 \\
          & \mathbf{1} & 1 & 0
        \end{array}
      \right]
     , \qquad 
     \widetilde{\IVEDual}
     = \left[
         \begin{array}{ccccc}
           \phantom{1}\\ 
           \mathbf{1} & 1 & 1 & 0 & \mathbf{1} \\
           1 & \mathbf{1} & 1 & 1 & 1 \\
           0 & 0 & 1 & 1 & 0
         \end{array}
       \right]
     ,
  \]
  with boldface entries representing the growth of the hyperedges and
  their duals through the elimination.
  Note that ``node elimination'' in a hypergraph is not the same as
  standard node elimination in a graph; it corresponds to merging
  the top row into all non-disjoint rows of the matrix $\Dual{\IVE}$.
\end{example}

\OldStuff{%
  Now consider "classical" node elimination as in the Gauss algorithm.
  We will analyze the relation between the fill-in evaluated in the
  node-node adjacency matrix of the associated graph and in the node-arc
  ncidence matrix of an associated hypergraph at a small example
  problem.
  
  \begin{example}
  Consider a problem represented by a graph $G=(V,E)$  with the node-node
  adjacency matrix $A^{\text{ad}}$ and the node-arc incidence matrix $A$
  as follows:
  $$A^{\text{ad}}=\left(\begin{array}{cccccc}
       1 & 1 & 0 & 0 & 0 & 1\\
       1 & 1 & 1 & 1 & 0 & 0\\
       0 & 1 & 1 & 0 & 1 & 0\\
       0 & 1 & 0 & 1 & 1 & 0\\
       0 & 0 & 1 & 1 & 1 & 1\\
       1 & 0 & 0 & 0 & 1 & 1
  \end{array}
  \right) \qquad 
  A=\left(\begin{array}{ccccccc}
       1 & 0 & 0 & 0 & 0 & 0 & 1\\
       1 & 1 & 1 & 0 & 0 & 0 & 0\\
       0 & 1 & 0 & 0 & 1 & 0 & 0\\
       0 & 0 & 1 & 1 & 0 & 0 & 0\\
       0 & 0 & 0 & 1 & 1 & 1 & 0\\
       0 & 0 & 0 & 0 & 0 & 1 & 1
  \end{array}
  \right)
  $$
  Suppose that we start with eliminating node $v_1$.
  $$A^{\text{ad}}=\left(\begin{array}{cccccc}
       \phantom{1}\\ 
        & 1 & 1 & 1 & 0 & \mathbf{1}\\
       \phantom{0} & 1 & 1 & 0 & 1 & 0\\
        & 1 & 0 & 1 & 1 & 0\\
        & 0 & 1 & 1 & 1 & 1\\
        & \mathbf{1} & 0 & 0 & 1 & 1
  \end{array}
  \right) \qquad 
  A=\left(\begin{array}{ccccccc}
       \phantom{1}\\ 
       1 & 1 & 1 & 0 & 0 & 0 & \mathbf{1}\\
       0 & 1 & 0 & 0 & 1 & 0 & 0\\
       0 & 0 & 1 & 1 & 0 & 0 & 0\\
       0 & 0 & 0 & 1 & 1 & 1 & 0\\
       \mathbf{1} & 0 & 0 & 0 & 0 & 1 & 1
  \end{array}
  \right)
  $$
  Since the hyperedges edge $e_1$ and $e_7$ in $A$ are identical, we
  will remove $e_7$ before the following step. In the next step, we
  eliminate $v_2$:
  $$A^{\text{ad}}=\left(\begin{array}{cccccc}
       \phantom{1}\\ 
       \phantom{1}\\ 
       \phantom{0} & \phantom{0} & 1 & \mathbf{1} & 1 & \mathbf{1}\\
        &  & \mathbf{1} & 1 & 1 & \mathbf{1}\\
        &  & 1 & 1 & 1 & 1\\
        &  & \mathbf{1} & \mathbf{1} & 1 & 1
  \end{array}
  \right) \qquad 
  A=\left(\begin{array}{ccccccc}
       \phantom{1}\\ 
       \phantom{1}\\ 
       \mathbf{1} & 1 & \mathbf{1} & 0 & 1 & 0 & \phantom{0}\\
       \mathbf{1} & \mathbf{1} & 1 & 1 & 0 & 0 & \\
       0 & 0 & 0 & 1 & 1 & 1 & \\
       {1} & \mathbf{1} & \mathbf{1} & 0 & 0 & 1 & 
  \end{array}
  \right)
  $$
  Before continuing with removing $v_3$, we again eliminate multiple
  hyperedges in $A$.
  $$A^{\text{ad}}=\left(\begin{array}{cccccc}
       \phantom{1}\\ 
       \phantom{1}\\ 
       \phantom{0}\\ 
        &  & \phantom{0} & 1 & 1 & {1}\\
        &  &  & 1 & 1 & 1\\
        &  &  & {1} & 1 & 1
  \end{array}
  \right) \qquad 
  A=\left(\begin{array}{ccccccc}
       \phantom{1}\\ 
       \phantom{1}\\ 
       \phantom{1}\\ 
       {1} & \phantom{0} & \phantom{0} & 1 & \mathbf{1} & 0 & \\
       \mathbf{1} &  & & 1 & 1 & 1 & \\
       {1} &  &  & 0 & \mathbf{1} & 1 & 
  \end{array}
  \right)
  $$
  \end{example}
  That does not seem to provide useful information.
}

There is another way to describe edge elimination in $G_{A}$ as node
elimination in a suitable graph, and since this corresponds to a square
matrix with symmetric nonzero pattern it allows to draw on the results
available for the solution of sparse symmetric positive definite linear
systems
  \cite{2006-Davis,1981-GeorgeLiu%
        }.
To this end we take a closer look at the edge--edge adjacency matrix
$\AEE$, more specifically at the process of running Gaussian elimination
on that matrix.


\section{Gaussian Elimination on the edge-edge adjacency matrix}%
\label{sec:edgeedge}

Let $G = (V,E)$ denote a hypergraph.
We now investgate how eliminating one of $G$'s edges changes the nonzero
pattern in the edge--edge adjacency matrix.

Let us first consider the symbolic elimination of an edge $x$, as
defined in
  \Cref{sec:seenhg}.
This elimination amounts to the following changes:
\[
  e \in E \setminus \{ x \}
  \rightarrow
  \begin{cases}
    e \cup x, & \text{if\ } e \cap x \neq \emptyset,\\
    e, & \text{else}.
  \end{cases}
\]
In particular this implies that all edges $e \in E \setminus\{x\}$ with
$e \cap x \neq \emptyset$ share all vertices of $x$ after its
elimination.
Thus, in terms of the edge--edge adjacency matrix $\AEE$, the
elimination results in a full block of nonzero entries covering all
$e \in E\setminus\{x\}$ with $e\cap x \neq \emptyset$.

On the other hand let us consider one step of symbolic Gaussian
elimination applied to the edge--edge adjacency matrix and note that
$\AEE$ is symmetric.
Without loss of generality let us assume that $\AEE$ is permuted such
that the edge $x$ is listed first.
Nonzero entries in the first column of $\AEE$ then correspond to edges
$e$ that share at least one vertex with $x$, i.e., for which
$e\cap x \neq \emptyset$.
Thus in the symbolic elimination step we now have to merge the nonzero
pattern of the first matrix row into the nonzero pattern of each row
corresponding to an edge $e$ with $e\cap x \neq \emptyset$.
Due to symmetry this again results in a full block of nonzeros
covering these edges (a clique in the graph $G_{\AEE}$ associated with
the matrix) and corresponds exactly to the nonzero pattern generated by
the symbolic edge elimination.

Thus in terms of the edge--edge connectivity structure, the symbolic
edge elimination process is equivalent to a symbolic Gaussian
elimination, applied to the edge--edge adjacency matrix.
Therefore this source of complexity, caused by increasing connectivity
among the remaining edges, can be approached in the same way it is done
in Gaussian elimination applied to sparse linear systems of equations.

Unfortunately, this does not cover all of the complexities of the
process.
If a fill-in element appears in $\AEE$ during Gaussian elimination
then this merely signals that all nodes from hyperedge $e_{j}$ will
be joined to those of $e_{i}$.
Therefore, the overall fill-in reflects the number of times when some
hyperedge will grow.
It does, however, not convey information about the current number of
nodes in the hyperedges, which would be necessary for assessing the cost
for the corresponding \ROneMod, see
  \cref{eq:costEdgeElimination}.


\section{NP-completeness results}
\label{sec:NP}

In this section we will show that even the problem
of minimizing the ``number of growths'' is NP-complete.

This follows directly from a well-known result stating the
NP-completeness of fill-in minimization
  \cite{1981-Yannakakis},
together with the following lemma.

\begin{lemma}
\label{lem:MatHypergraph}

  The nonzero pattern of any symmetric positive definite irreducible
  $n$-by-$n$ matrix can be interpreted as the edge--edge adjacency
  matrix of a suitable hypergraph $G = ( V, E )$ with $| E | = n$
  edges.

  \begin{proof}
    Define
    \[
      V = \{ v_{i,j} \mid i > j, a_{i,j} \not= 0 \}
      ,
    \]
    that is, we have one node for each nonzero in the strict lower
    triangle of $A$.
    Let $E = \{ e_{1}, \ldots, e_{n} \}$, where
    \[
      e_{j} = \{ v_{i,j} \mid i > j, a_{i,j} \not= 0 \}
              \cup \{ v_{j,i} \mid j > i, a_{j,i} \not= 0 \}
      ,
    \]
    i.e., $e_{j}$ contains just those nodes corresponding to nonzeros
    in column $j$ or row $j$ of $A$'s strict lower triangle.
    Note that $e_{j} \not= \emptyset$ because otherwise row and
    column $j$ of $A$ would contain just the diagonal entry, i.e.,
    $A$ were reducible.
    Then, for $k > j$ we have
    \[
        e_{k} \cap e_{j}
      =
        \{ v_{i,j} \mid i > j, a_{i,j} \not= 0 \}
        \cap
        \{ v_{k,i} \mid k > i, a_{k,i} \not= 0 \}
    \]
    (the other three intersections being empty), and this is nonempty
    iff there is a node $v_{i,j} \equiv v_{k,i}$ in both column $j$
    and row $k$, i.e., $a_{k,j} \not= 0$.
    Using Lemma~\ref{lem:incidence-adjacency},
    this implies that $A$ and $\AEE = \Herm{\IVE} \IVE$ have the same
    nonzero pattern.
  \end{proof}
\end{lemma}

\begin{remark}
  In most cases, the same nonzero pattern may also be obtained with
  hypergraphs containing fewer nodes.
  It is therefore tempting to take $\IVE$ to be the nonzero pattern of
  the Cholesky factor $U$ from $A = \Herm{U} U$ in order to obtain the
  sparsity pattern of $A$ with a hypergraph containing just $n$ nodes.
  Unfortunately, cancellation in the product $\Herm{U} U$ may introduce
  zeros in $A$ that are not present in the product $\Herm{\IVE} \IVE$
  obtained this way, and this cancellation can be structural.
  In fact, exhaustive search reveals that, for $n = 5$, the pattern
  \[
    A = \left[
          \begin{array}{ccccc}
            1 & 0 & 1 & 1 & 1 \\
            0 & 1 & 1 & 1 & 1 \\
            1 & 1 & 1 & 0 & 0 \\
            1 & 1 & 0 & 1 & 0 \\
            1 & 1 & 0 & 0 & 1
          \end{array}
        \right]
  \]
  cannot be obtained as $\Herm{\IVE} \IVE$ with any hypergraph
  containing fewer than six nodes, and six nodes are sufficient
  according to the proof of Lemma~\ref{lem:MatHypergraph}
  because the strict lower triangle of $A$ contains six nonzeros.
\end{remark}

Note that for Lemma~\ref{lem:MatHypergraph}
we have assumed that we may start with a hypergraph; cf.\ Remark~\ref{rem:MoreGeneralSettings}.
If this is not allowed and we restrict ourselves to eliminating ``true
edges,'' thus zeroing one pair of matrix entries $a_{k,\ell}$ and
$a_{\ell,k}$ per step, then a simple combinatorial argument shows that
there must be symmetric positive definite matrices whose nonzero pattern
cannot be interpreted as that of an edge--edge adjacency matrix
$\AEE = \IVE^{T} \cdot \IVE$ to \emph{any} graph $G = ( V, E )$.

To see this, we note that the number of nonzero patterns for a symmetric
$n$-by-$n$ matrix $A$ is
$\nu_{A} = 2^{n(n-1)/2} = \left( 2^{(n-1)/2} \right)^{n}$,
because each of the $n(n-1)/2$ entries in the strict lower triangle may
be zero or not.
Now assume that the matrix has the same nonzero pattern as
$\IVE^{T} \cdot \IVE$ for some graph $G = ( V, E )$ with $n$ edges and
some number of nodes, $v$.
Then $\IVE \in \MathR^{v \times n}$ contains exactly two nonzeros in
each of its columns, and we may assume w.l.o.g.\ that $v \leq 2 n$,
because at most $2 n$ rows of $\IVE$ can contain a nonzero, and rows
with all zeros can be removed without affecting the product
$\IVE^{T} \cdot \IVE$ (this corresponds to removing isolated nodes from
$G$).
Then there are at most ${2n \choose 2} = 2n ( 2n - 1 ) / 2$ possible
combinations for the positions of the two nonzeros in each column of
$\IVE$, leading to the overall number of possible matrices $\IVE$ being
bounded by $\nu_{\IVE} \leq \left( \frac{2n ( 2n - 1 )}{2} \right)^{n}$.
Since $2^{(n-1)/2} > \frac{2n ( 2n - 1 )}{2}$ for large $n$,
we also have $\nu_{A} > \nu_{\IVE}$, and therefore not all symmetric
matrices can be interpreted as edge--edge adjacency matrices.

In this situation the proof of NP-completeness for fill-in minimization
does not carry over, and it is currently not known whether this
restricted problem is indeed NP-complete.

In the light of these results one still may try to find orderings that
lead to reduced (arithmetic or memory) complexity without being optimal
in the above sense.
This will be discussed in the following.


\section{Heuristics for choosing edge elimination orderings}
\label{sec:seecomp}

Based on the findings in
  \Cref{sec:see,sec:seenhg}
it is natural to analyze the complexity of Algorithm~\ref{alg:edgeElimination}
in terms of the overall number of roots of the secular equation that
have to be calculated during all edge eliminations.
Combining this analysis with the cost for the calculation of a single
root of the secular equation gives us direct access to the complexity
of the \Hermit{}ian (sparse) eigenvalue problem.

\begin{lemma}
  Let $G_{A} = (V,E)$ be an undirected graph of a matrix $A$,
  interpreted as a hypergraph.
  Further define an ordering of the edges $e_{1},\ldots,e_{|E|}$.
  Then the total number $N$ of secular equation roots that have to be
  calculated in Algorithm~\ref{alg:edgeElimination}
  is given by
  \[
    N = \sum_{j = 1}^{| E |} N_{e_{j}}
    ,
  \]
  using the definition of $N_{e_{j}}$ from
    \cref{eq:costEdgeElimination}.
\end{lemma}

\paragraph{Minimum incidence (MI) ordering}

In analogy to the minimum degree ordering in Gaussian elimination the
first heuristic that comes to mind accounts for the number of incident
edges. In the hypergraph setting two edges $e$ and $x$ are incident iff
$e \cap x \neq \emptyset$, i.e., when eliminating $x$ the edge $e$
changes and vice versa.
By introducing the quantities
\[
  \mu_{\rm i}(x)
  =
  | \{e \in E \mid e \cap x \neq \emptyset \} |
  ,
\]
the strategy thus chooses in every step the edge with the fewest
incident edges.
Once an edge $x$ is eliminated, the number of incident edges needs to be
updated only for all edges $e$ that have been incident with $x$.

\paragraph{Minimal root number (MR) ordering}

Another heuristic is to account for the number of roots of the secular equation that need to be calculated when eliminating a hyperedge. That is, we define the quantities
\[
  \mu_{\rm r}(x) = |x|
  ,
\]
and the MR strategy chooses in every step the edge with the smallest number of contained vertices. After elimination of an edge, $\mu_{\rm r}$ needs to be updated for all edges incident with the eliminated edge.

\paragraph{Minimal roots/costs with look-ahead (MC)}

The last heuristic under consideration modifies the MR heuristic by adding a look-ahead component. The elimination of an edge $x$ incurs a growth of all edges $e$ with $x\cap e \neq \emptyset$ by $|x \cup e| - |e|$ vertices. This in turn relates to the number of roots that need to be calculated in a future elimination. Due to the fact that the cost of eliminating an edge $x$ with $|x|$ nodes is proportional to $|x|^2$ we consider the two measures
\[
\mu_{\rm c}^{(k)}(x) = |x|^k + \sum_{e\cap x \neq \emptyset} |x \cup e|^k - |e|^k
\] for $k = 1,2$ and choose to eliminate the edge with the current smallest value of $\mu_{\rm c}^{(k)}$. Due to the look-ahead nature of the measure, updating it now involves not only the edges incident with $x$, but also the next-neighbors as well.

In order to assess the efficiency of these heuristics, they have been
applied to matrices with different sparsity patterns, i.e., different
structures of the associated graph $G_{A}$.


\paragraph{I. The chain graph}
In order to enable a comparison of our approach to the tridiagonal divide-and-conquer algorithm we first apply the symbolic process to a chain of $N$ nodes, which is the graph corresponding to a tridiagonal matrix. The divide-and-conquer strategy for this graph results in the calculation $N$ roots on each level of the recursion for a total of $N\log_2(N)$ roots.

As can be seen from the results in Table~\ref{tab:chain_results} both the strategy that chooses the edge with currently smallest number of contained vertices, based on $\mu_{\rm r}$, as well as the strategy that accounts for the current and future cost of eliminating an edge, based on $\mu_{\rm c}^{(1)}$, result in elimination orderings which are equivalent to the divide-and-conquer strategy. While the strategy based on measure $\mu_{\rm c}^{(2)}$ comes close to the optimal total number of roots, the strategy based on choosing to eliminate the edge with the least number of incident edges fails spectacularly and eliminates the edges in lexicographic ordering.

\begin{table}[]
    \caption{Symbolic elimination for the chain graph with $N = 256$ nodes. Reported are the accumulated number of roots that need to be calculated over the whole elimination process. For comparison, the number of root calculations in the divide and conquer algorithm for this problem is $256\times \log_2(256) = 2048$.}
        \label{tab:chain_results}
    \centering
    \begin{tabular}{c|c|c|c|c}
    Heuristic & $\mu_{\rm i}$ & $\mu_{\rm r}$ & $\mu_{\rm c}^{(1)}$ & $\mu_{\rm c}^{(2)}$ \\\hline
    $\sum |x|$ & $16766$ & $2048$ & $2048$ & $2152$
    \end{tabular}
    
\end{table}

\begin{figure}[!h]
  \begin{center}
    \includegraphics[width=\textwidth]{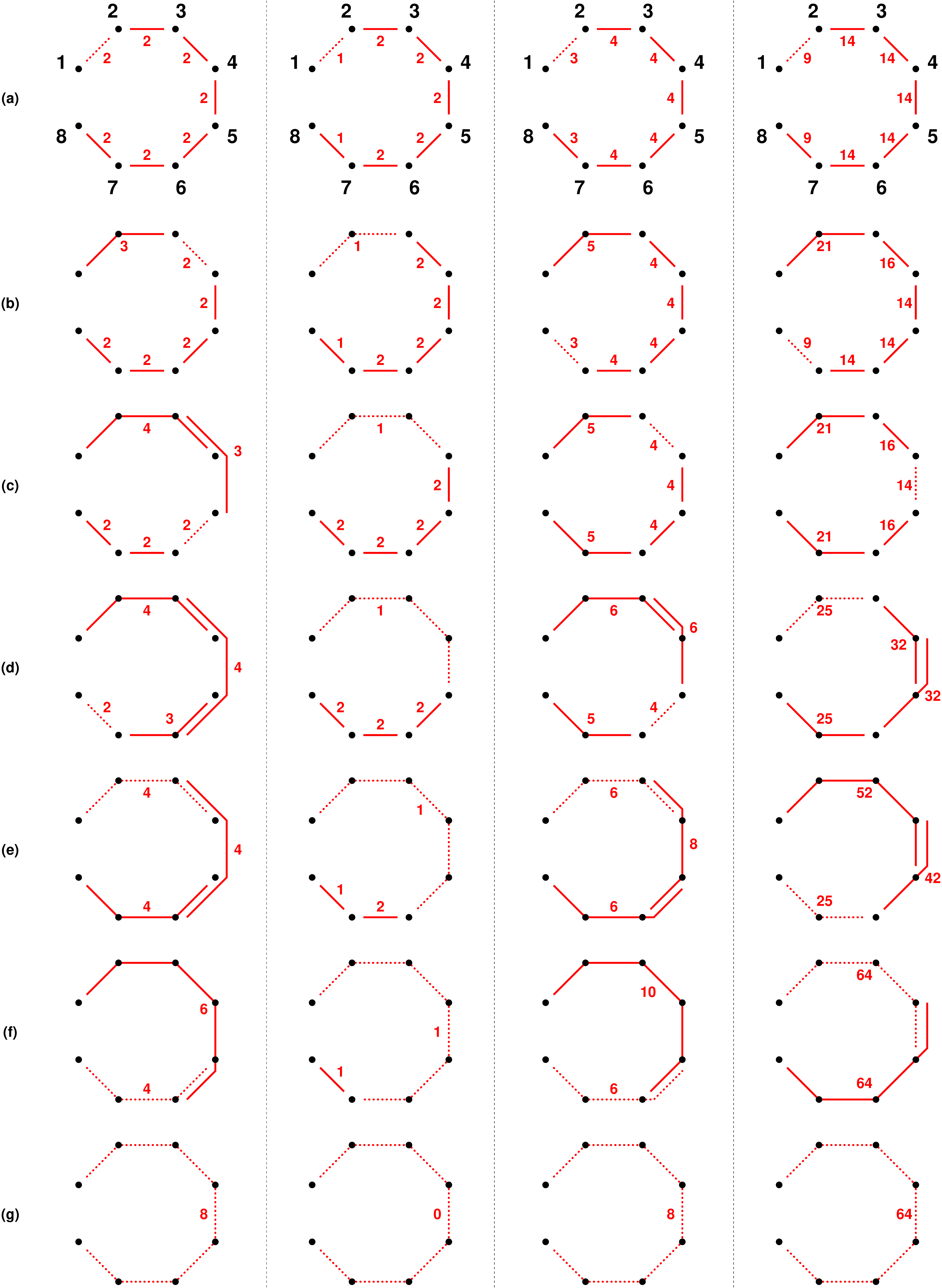}
  \end{center}
  \caption{Order of (hyper)edge elimination for a size-$8$ tridiagonal
           matrix with the strategies (from left to right)
           Minimal root number (MR,
             minimize $\mu_{\rm r}( x )$),
           Minimum incidence (MI,
             minimize $\mu_{\rm i}( x )$),
           Minimal roots with look-ahead (MC1,
             minimize $\mu_{\rm c}^{(1)}( x )$), and
           Minimal costs with look-ahead (MC2,
             minimize $\mu_{\rm c}^{(2)}( x )$).
           For each of the seven elimination steps (a) to (g),
           the remaining (hyper)edges are shown together with their
           $\mu$ values, and the (hyper)edge selected for elimination is
           highlighted as a dotted line.
           If the minimum is not unique then the ``first'' minimizing
           hyperedge (clockwise) is chosen for elimination.}
  \label{fig:HeuristicsChain8}
\end{figure}

The progress of the elimination for a chain graph with $N = 8$ nodes
is shown in Figure~\ref{fig:HeuristicsChain8}.
Again, $\mu_{\rm r}$ and $\mu_{\rm c}^{(1)}$ achieve
the same $\sum |x|$ value as tridiagonal divide-and-conquer,
$N \log_{2}( N ) = 24$, $\mu_{\rm c}^{(2)}$ is slightly
worse ($\sum |x| = 25$), and $\mu_{\rm i}$ leads to the
lexicographic ordering ($\sum |x| = 35$).

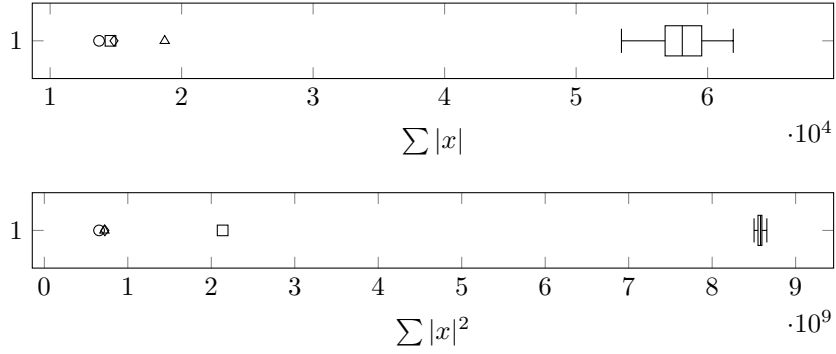
\begin{figure}[!h]
\begin{center}
    \begin{tikzpicture}
        \pgfplotstableread{stats_N256_lattice_baseline.dat}{\data}
        \pgfplotstableread{stats_N256_lattice_heuristics.dat}{\dataH}
        \begin{axis}[
	    	boxplot/draw direction = x,
		    xlabel={$\sum|x|$},
		    ytick={1},
		    y=.5cm,
		    width=\linewidth,
		    ymin=0,
		    ymax=2,
        	]
        	\addplot[boxplot] table[y=linear] {\data};
	        \pgfplotstablegetelem{0}{linear}\of\dataH   
	        \pgfmathsetmacro\tmpval{\pgfplotsretval}
	        \addplot[mark=triangle] coordinates {(\tmpval,1)};
	        \pgfplotstablegetelem{1}{linear}\of\dataH   
	        \pgfmathsetmacro\tmpval{\pgfplotsretval}
	        \addplot[mark=square] coordinates {(\tmpval,1)};
	        \pgfplotstablegetelem{2}{linear}\of\dataH   
	        \pgfmathsetmacro\tmpval{\pgfplotsretval}
	        \addplot[mark=o] coordinates {(\tmpval,1)};
	        \pgfplotstablegetelem{3}{linear}\of\dataH   
	        \pgfmathsetmacro\tmpval{\pgfplotsretval}
	        \addplot[mark=diamond] coordinates {(\tmpval,1)};
	    \end{axis}
    \end{tikzpicture}\vspace{1em}
    \begin{tikzpicture}
        \pgfplotstableread{stats_FEM_baseline.dat}{\data}
        \pgfplotstableread{stats_FEM_heuristics.dat}{\dataH}
        \begin{axis}[
	    	boxplot/draw direction = x,
	        xlabel={$\sum|x|^2$},
		    ytick={1},
		    y=.5cm,
		    width=\linewidth,
		    ymin=0,
		    ymax=2,
        	]
        	\addplot[boxplot] table[y=quadratic] {\data};
	        \pgfplotstablegetelem{0}{quadratic}\of\dataH   
	        \pgfmathsetmacro\tmpval{\pgfplotsretval}
	        \addplot[mark=triangle] coordinates {(\tmpval,1)};
	        \pgfplotstablegetelem{1}{quadratic}\of\dataH   
	        \pgfmathsetmacro\tmpval{\pgfplotsretval}
	        \addplot[mark=square] coordinates {(\tmpval,1)};
	        \pgfplotstablegetelem{2}{quadratic}\of\dataH   
	        \pgfmathsetmacro\tmpval{\pgfplotsretval}
	        \addplot[mark=o] coordinates {(\tmpval,1)};
	        \pgfplotstablegetelem{3}{quadratic}\of\dataH   
	        \pgfmathsetmacro\tmpval{\pgfplotsretval}
	        \addplot[mark=diamond] coordinates {(\tmpval,1)};
	    \end{axis}
    \end{tikzpicture}
    \end{center}
    \caption{Accumulated number of roots $\sum |x|$ (top) and root calculation costs $\sum |x|^2$ (bottom) for a regular $16\times 16$ lattice graph with $N = 256$ nodes.
    Results for the heuristics are plotted as $(\mu_{\rm i},\triangle)$, $(\mu_{\rm r},\square)$, $(\mu_{\rm c}^{(1)},\circ)$ and $(\mu_{\rm c}^{(2)},\diamond)$ (towards the left),
    and the boxplots close to the right summarize the results for 20 random elimination orderings.}
    \label{fig:lattice_N256}
\end{figure}

\paragraph{II. Structured graphs} 
Structured graphs are often encountered in discretizations of partial differential equations. The resulting graphs are planar and usually possess a large diameter. In Figure~\ref{fig:lattice_N256} we report results in terms of accumulated number of roots $\sum |x|$ and cost of root elimination $\sum|x|^2$ of the hypergraph edge elimination approach for a uniform $16\times 16$ lattice. We compare the results for the four heuristics with a statistical baseline of $20$ random elimination orderings. As can be seen from the figure all four heuristics yield largely reduced cost measures compared to the baseline. Notably, the ordering of the heuristics in terms of the two cost measures is not identical, i.e., an overall minimal number of accumulated root calculations does not immediately lead to a minimal accumulated root elimination cost.

\begin{figure}[!h]
    \begin{center}
    \begin{tikzpicture}
        \pgfplotstableread{stats_FEM_baseline.dat}{\data}
        \pgfplotstableread{stats_FEM_heuristics.dat}{\dataH}
        \begin{axis}[
	    	boxplot/draw direction = x,
		    xlabel={$\sum|x|$},
		    ytick={1},
		    y=.5cm,
		    width=\linewidth,
		    ymin=0,
		    ymax=2,
        	]
        	\addplot[boxplot] table[y=linear] {\data};
	        \pgfplotstablegetelem{0}{linear}\of\dataH   
	        \pgfmathsetmacro\tmpval{\pgfplotsretval}
	        \addplot[mark=triangle] coordinates {(\tmpval,1)};
	        \pgfplotstablegetelem{1}{linear}\of\dataH   
	        \pgfmathsetmacro\tmpval{\pgfplotsretval}
	        \addplot[mark=square] coordinates {(\tmpval,1)};
	        \pgfplotstablegetelem{2}{linear}\of\dataH   
	        \pgfmathsetmacro\tmpval{\pgfplotsretval}
	        \addplot[mark=o] coordinates {(\tmpval,1)};
	        \pgfplotstablegetelem{3}{linear}\of\dataH   
	        \pgfmathsetmacro\tmpval{\pgfplotsretval}
	        \addplot[mark=diamond] coordinates {(\tmpval,1)};
	    \end{axis}
    \end{tikzpicture}\vspace{1em}
    \begin{tikzpicture}
        \pgfplotstableread{stats_FEM_baseline.dat}{\data}
        \pgfplotstableread{stats_FEM_heuristics.dat}{\dataH}
        \begin{axis}[
	    	boxplot/draw direction = x,
	        xlabel={$\sum|x|^2$},
		    ytick={1},
		    y=.5cm,
		    width=\linewidth,
		    ymin=0,
		    ymax=2,
        	]
        	\addplot[boxplot] table[y=quadratic] {\data};
	        \pgfplotstablegetelem{0}{quadratic}\of\dataH   
	        \pgfmathsetmacro\tmpval{\pgfplotsretval}
	        \addplot[mark=triangle] coordinates {(\tmpval,1)};
	        \pgfplotstablegetelem{1}{quadratic}\of\dataH   
	        \pgfmathsetmacro\tmpval{\pgfplotsretval}
	        \addplot[mark=square] coordinates {(\tmpval,1)};
	        \pgfplotstablegetelem{2}{quadratic}\of\dataH   
	        \pgfmathsetmacro\tmpval{\pgfplotsretval}
	        \addplot[mark=o] coordinates {(\tmpval,1)};
	        \pgfplotstablegetelem{3}{quadratic}\of\dataH   
	        \pgfmathsetmacro\tmpval{\pgfplotsretval}
	        \addplot[mark=diamond] coordinates {(\tmpval,1)};
	    \end{axis}
    \end{tikzpicture}
    \end{center}
    \caption{Accumulated number of roots $\sum |x|$ (top) and root calculation costs $\sum |x|^2$ (bottom) for a triangulation of the unit disc with $N = 1313$ nodes. Results for the heuristics are plotted as $(\mu_{\rm i},\triangle)$, $(\mu_{\rm r},\square)$, $(\mu_{\rm c}^{(1)},\circ)$ and $(\mu_{\rm c}^{(2)},\diamond)$. The boxplots represent results of $20$ random elimination orderings.}
    \label{fig:triangulation_random}
\end{figure}
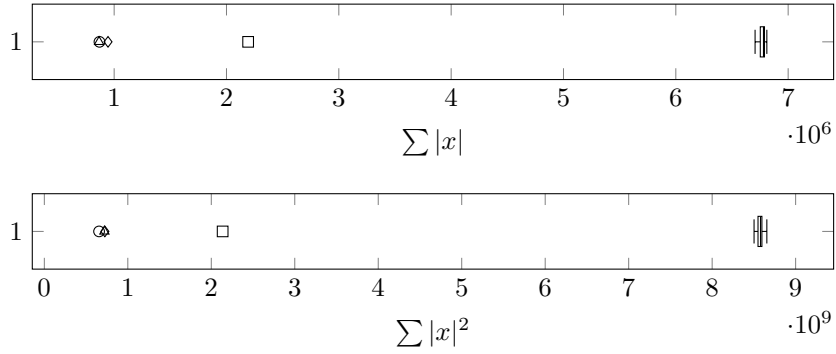

Next we apply the same test setup to a graph that is a triangulation of the unit disc with $1313$ nodes. In Figure~\ref{fig:triangulation_random} we report accumulated number of roots $\sum |x|$ and root elimination costs $\sum|x|^2$ for the four heuristics and report the statistical baseline of $20$ random orderings. Again we see that all four heuristics are clearly better than using a random elimination ordering.

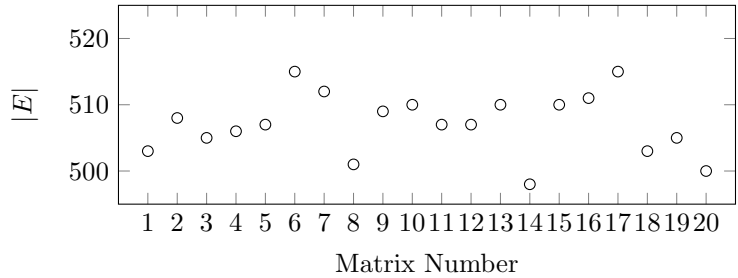
\begin{figure}[!h]
    \begin{center}
    \begin{tikzpicture}
    \pgfplotstableread{stats_N128_nnz08_matrix.dat}{\data}
    \begin{axis}[
	    xlabel={Matrix Number},
	    ylabel={$|E|$},
		xtick={1,...,20},
		ytick={500,510,520},
		width=.8\linewidth,
		height=12em,
		xmin=0,
		xmax=21,
		ymin=495,
		ymax=525
        ]
        \foreach \i in {0,1,...,19}{
	            \pgfplotstablegetelem{\i}{edges}\of\data   
	            \pgfmathsetmacro\tmpval{\pgfplotsretval}
                \addplot[mark=o] coordinates {({\i+1},\tmpval)};
	        }
	        
	    \end{axis}
    \end{tikzpicture}
    \end{center}
    \caption{Number of edges $|E|$ of $20$ randomly generated sparse graphs.}
    \label{fig:number_edges_random_graphs}
\end{figure}

\paragraph{III. Sparse random graphs} Finally we compare the heuristics for randomly generated graphs. We use the Matlab built-in function \texttt{sprandsym} to generate an undirected graph with $N$ nodes with a non-zero density of $\tfrac{8}{N}$. The resulting graphs' average degree is thus approximately $8$. We now test the heuristics for $20$ such graphs of sizes $N = 128$. 
In Figure~\ref{fig:number_edges_random_graphs} we report the number of edges of the matrices used in the tests.

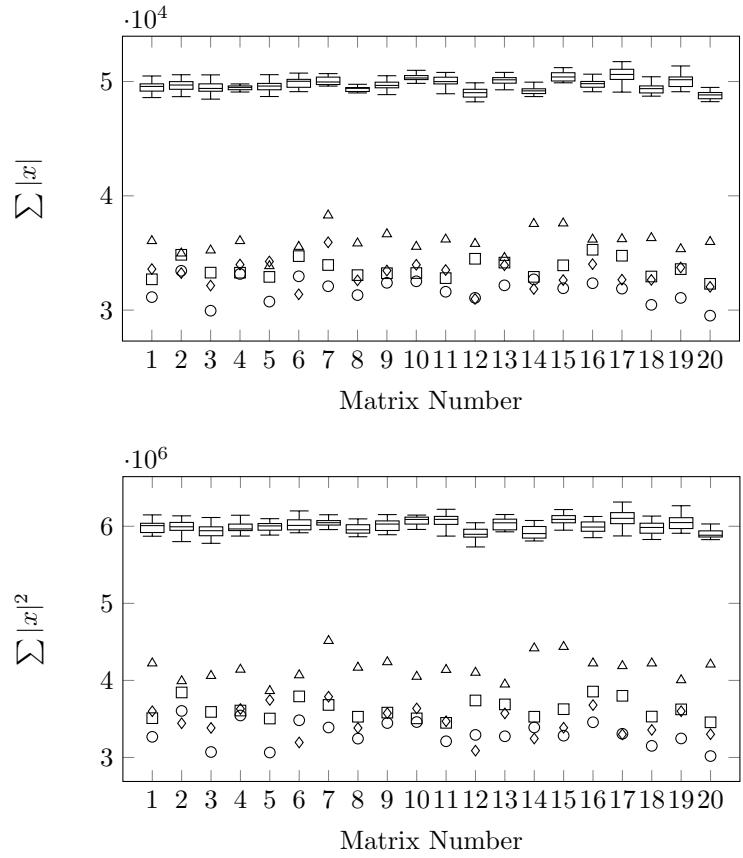
\begin{figure}
\begin{center}
\begin{tikzpicture}
    \pgfplotstableread{stats_N128_nnz08_baseline.dat}{\data}
    \pgfplotstableread{stats_N128_nnz08_heuristics.dat}{\dataH}
        \begin{axis}[
	    	boxplot/draw direction = y,
		    xlabel={Matrix Number},
		    ylabel={$\sum|x|$},
		    xtick={1,...,20},
		    width=.8\linewidth,
		    height=16em,
		    xmin=0,
		    xmax=21,
        	]
        	\foreach \i in {1,2,...,20}{
	            \addplot[boxplot] table[y=l\i] {\data};
	            \pgfplotstablegetelem{0}{l\i}\of\dataH   
	            \pgfmathsetmacro\tmpval{\pgfplotsretval}
	            \addplot[mark=triangle] coordinates {(\i,\tmpval)};
	            \pgfplotstablegetelem{1}{l\i}\of\dataH   
	            \pgfmathsetmacro\tmpval{\pgfplotsretval}
	            \addplot[mark=square] coordinates {(\i,\tmpval)};
	            \pgfplotstablegetelem{2}{l\i}\of\dataH   
	            \pgfmathsetmacro\tmpval{\pgfplotsretval}
	            \addplot[mark=o] coordinates {(\i,\tmpval)};
	            \pgfplotstablegetelem{3}{l\i}\of\dataH   
	            \pgfmathsetmacro\tmpval{\pgfplotsretval}
	            \addplot[mark=diamond] coordinates {(\i,\tmpval)};
	        }
	    \end{axis}
    \end{tikzpicture}\vspace{1em}
    \begin{tikzpicture}
    \pgfplotstableread{stats_N128_nnz08_baseline.dat}{\data}
    \pgfplotstableread{stats_N128_nnz08_heuristics.dat}{\dataH}
        \begin{axis}[
	    	boxplot/draw direction = y,
		    xlabel={Matrix Number},
		    ylabel={$\sum|x|^2$},
		    xtick={1,...,20},
		    width=.8\linewidth,
		    height=16em,
		    xmin=0,
		    xmax=21,
        	]
        	\foreach \i in {1,2,...,20}{
	            \addplot[boxplot] table[y=q\i] {\data};
	            \pgfplotstablegetelem{0}{q\i}\of\dataH   
	            \pgfmathsetmacro\tmpval{\pgfplotsretval}
	            \addplot[mark=triangle] coordinates {(\i,\tmpval)};
	            \pgfplotstablegetelem{1}{q\i}\of\dataH   
	            \pgfmathsetmacro\tmpval{\pgfplotsretval}
	            \addplot[mark=square] coordinates {(\i,\tmpval)};
	            \pgfplotstablegetelem{2}{q\i}\of\dataH   
	            \pgfmathsetmacro\tmpval{\pgfplotsretval}
	            \addplot[mark=o] coordinates {(\i,\tmpval)};
	            \pgfplotstablegetelem{3}{q\i}\of\dataH   
	            \pgfmathsetmacro\tmpval{\pgfplotsretval}
	            \addplot[mark=diamond] coordinates {(\i,\tmpval)};
	        }
	        
	    \end{axis}
    \end{tikzpicture}
    \end{center}
\caption{Accumulated number of roots $\sum |x|$ (top) and root calculation costs $\sum |x|^2$ (bottom) for the $20$ randomly generated sparse graphs. Results for the heuristics are plotted as $(\mu_{\rm i},\triangle)$, $(\mu_{\rm r},\square)$, $(\mu_{\rm c}^{(1)},\circ)$ and $(\mu_{\rm c}^{(2)},\diamond)$. Each boxplot represents results of $20$ random elimination orderings.}
    \label{fig:results_random_graphs_k8}
\end{figure}

In Figure~\ref{fig:results_random_graphs_k8} we report the results of the heuristics applied to these randomly generated sparse graphs. We report both the accumulated number of roots $\sum |x|$ as well as the accumulated cost of root calculations $\sum |x|^2$. In order to gauge the potential gains realized by the heuristics, we include boxplots of $20$ random elimination orderings as well.

Overall, our experiments suggest that, while none of the proposed
strategies is consistently superior, choosing the hyperedge with
minimum $\mu_{\rm c}^{(1)}$ value for elimination seems to be a
reasonable way to reduce both cost measures, the total number of roots
to compute, $\sum | x |$, and the operations to do this,
$\sum | x |^{2}$.


\section{Concluding remarks}\textbf{}
\label{sec:conclusions}
We have shown in this paper that the symmetric eigenvalue problem can be interpreted as an elimination process, where all edges of the corresponding graph need to be eliminated. This symbolic equivalence is facilitated by a hypergraph point of view and in complete analogy to the vertex elimination that characterizes the symbolic solution of linear systems by means of Gaussian elimination.

Furthermore, we showed that the hypergraph information in every stage of the elimination process is captured by symbolic Gaussian elimination applied to the edge--edge adjacency matrix---a formal dual to the regular vertex--vertex adjacency matrix. Exploiting this connection we were able to transfer the result of NP-hardness for the calculation of an optimal elimination ordering from the linear systems case to the symmetric eigenvalue problem.

While optimality cannot be achieved, we proposed different heuristics to determine good elimination orderings and numerically explored their use. In particular, we compared them to a baseline of random elimination orderings, where they proved to be vastly superior to this baseline. We also explored if the chosen heuristics are able to reproduce the optimal ordering in case that the graph of the matrix is a chain graph, i.e., the matrix is tri-diagonal. In this case, the proposed edge elimination algorithm with optimal elimination ordering is equivalent to an iterative (rather than recursive) formulation of the divide-and-conquer approach to tridiagonal symmetric eigenvalue problems.

Considered from the point of view of this paper, the usual approach of initial reduction to tridiagonal form and subsequent solution of the tridiagonal eigenvalue problem can be viewed as the reduction to a chain graph with subsequent edge elimination, for which an optimal elimination strategy is known.

The equivalence of the Hermitian eigenvalue problem and symbolic hypergraph edge elimination can be easily transferred to the calculation of the singular value decompostion based on the observation the the singular value decompostion $AV = U\Sigma$ of $A \in \mathbb{C}^{m\times n}$ can be computed by considering the Hermitian eigenvalue problem
\[
\begin{bmatrix} 0 & A^{H}\\A & 0\end{bmatrix} \begin{bmatrix} V & V\\ U & -U\end{bmatrix} = \begin{bmatrix} V & V\\ U & -U\end{bmatrix} \begin{bmatrix} \Sigma & 0\\0 & -\Sigma \end{bmatrix}.
\]

\section*{Acknowledgements} The authors would like to thank Kathrin Klamroth and Michael Stiglmayr for many fruitful discussions.



\bibliographystyle{siam}
\bibliography{references}
\end{document}